\documentclass[11pt]{article}
\usepackage{amsmath,amsthm,amsfonts,amssymb,amscd}
\usepackage[latin1]{inputenc}

\usepackage{graphicx}

\usepackage{amsmath}

\usepackage{amssymb}

\usepackage{amscd}

\usepackage{bbm}

\newtheorem{Theorem}{Theorem}
\newtheorem{Lemma}{Lemma}

\newtheorem{Proposition}{Proposition}

\newtheorem{Problem}{Problem}
\newtheorem{prope}{Property}

\newtheorem{Definition}{Definition}
\newenvironment{definition}{\begin{Definition} \em}{\end{Definition}}

\newtheorem{Remark}{Remark}
\def\O{{\mathcal{O}}}

\def\ov{\overline}

\def\bc{{\mathbb{C}}}

\def\Gen{\operatorname{{Gen}}}

\def\Hol{\operatorname{{Hol}}}

\def\Ker{\operatorname{{Ker}}}

\def\Diff{\operatorname{{Diff}}}
\def\sing{\operatorname{{sing}}}
\def\sing{\operatorname{{sing}}}
\def\codim{\operatorname{{codim}}}

\def\Sat{\operatorname{{Sat}}}

\begin{document}

\title{On the  integrability of holomorphic vector fields}
\author{L. C\^{a}mara and B. Sc\'{a}rdua}
\date{}
\maketitle

\begin{abstract}
We determine topological and algebraic conditions for a germ of
holomorphic foliation $\mathcal F(X)$ induced by a generic vector
field $X$ on $(\mathbb{C}^{3},0)$ to have a holomorphic first
integral, i.e., a germ of holomorphic map $F
\colon(\mathbb{C}^{3},0)\longrightarrow(\mathbb{C}^{2},0)$ such
that the leaves of $\mathcal F(X)$ are contained in the level
curves of $F$.

\end{abstract}

%\tableofcontents

\section{Introduction\label{section:introduction}}

This paper is devoted to the relations between  certain algebraic,
topological and geometric properties of holomorphic
foliations\footnote{2000 \textit{Mathematics Subject
Classification}. Primary 37F75; Secondary 32M25, 32S65.}. The
question of existence of first integrals for holomorphic
foliations on $(\mathbb{C}^{2},0)$ depends heavily on the
resolution algorithm due to Bendixson (\cite{Be 01}) and after
rediscovered by Seidenberg (\cite{Se 68}). In dimension three the
subject has been developed by Cano (\cite{Ca 87}) and more
recently, in the real analytic case, by Panazzolo (\cite{Pa 05}),
with the aid of the concept of weighted blow-up. This resolution
leads to an ambient three-fold which is no more regular but has
toroidal singularities (more precisely weighted projective
spaces). On the other hand, as a consequence of a monomialization
theorem due to Cutkosky (\cite{Cuk 2002}), any singular
holomorphic foliation on $(\mathbb{C}^{3},0)$ having a holomorphic
first integral (see definition below) can be resolved with a
finite number of monoidal transformations, avoiding the presence
of singularities in the ambient three-fold of the resolved
foliation. Generically, the singularities which appear after the
resolution process are non-degenerate (and thus isolated), and
have three (transverse) invariant planes through each singularity.
Due to this fact, we consider along this work mainly singular
holomorphic foliations satisfying the above properties.

In this work we study topological and algebraic conditions for the
existence of holomorphic first integrals in a spirit similar to
 to the classical approach for the two-dimensional case
studied  in \cite{MaMo 80}. Neverteless, some differences arise
with respect to the $2$-dimensional case. The first is the fact
that a stronger condition than the Siegel condition on the
eigenvalues is required (see
Definition~\ref{definition:conditionstar} below). Also in order to
find suitable integrable codimension one distributions containing
the trajectories of $X$  we need to apply some techniques from
Partial Differential Equations (Problem~\ref{geometric
integrability} and \S \ref{subsection:solution}).
\subsection{Notation and statements}
\label{subsection:notation}
Denote the ring of germs of holomorphic functions on $(\mathbb{C}%
^{n},0)$ by $\mathcal{O}_{n}$ and its maximal ideal by $\mathcal{M}_{n}$. We
 say that a germ of one-dimensional holomorphic foliation $\mathcal F(X)$
 on $(\mathbb{C}^{n},0)$ given by the (germ of) holomorphic
vector field $X\in\mathcal{X}(\mathbb{C}^{n},0)$ is
\textit{integrable}, or  \textit{has a holomorphic first
integral}, if there is a germ of holomorphic map
$F:(\mathbb{C}^{n},0)\longrightarrow(\mathbb{C}^{n-1},0)$ such
that $(a)$ the leaves of $\mathcal F(X)$ are contained in level
curves of $F$ or equivalently $i_{X}dF=0$ and $(b)$ $F$ is a
submersion off $\sing(\mathcal F(X))$,
the \textit{singular set of }$\mathcal F(X)$. Further, let $f:(\mathbb{C}%
^{n},0)\longrightarrow(\overline{\mathbb{C}},\infty)$ be a
meromorphic function, then we  say that $f$ is $\mathcal
F(X)$\textit{-invariant} if $i_{X}df=0$, i.e.,if the leaves of
$\mathcal F(X)$ are contained in the level sets of $f$. A germ of
vector field on $(\mathbb{C}^{n},0)$ is \textit{non-degenerate} if
$DX(0)$ has just non-vanishing eigenvalues. Recall that
generically $DX(0)$ has three distinct eigenvalues and thus is
diagonalizable and $X$ has an isolated singularity at the origin.
Furthermore, from Poincar\'{e}-Dulac, Siegel and Brjuno
linearization theorems and from \cite{CaKuiPa 78}, generically $X$
leaves invariant the \textit{coordinate hyperplanes} $x_{1}\cdots
x_{n}=0$.. Therefore we  say that $\mathcal F(X)$ is
\textit{non-degenerate generic} if $DX(0)$ is diagonalizable and
$X$ leaves invariant the coordinate planes. For simplicity, we
denote the set of germs of non-degenerate generic vector fields on
$(\mathbb{C}^{n},0)$ by $\Gen(\frak X(\bc^n,0))$. Let
$X\in\Gen(\frak X(\bc^n,0))$, $S$ a smooth integral curve of
$\mathcal F(X)$ through the origin, and $f$ a germ of $\mathcal
F(X)$-invariant meromorphic function. Then we denote by
$\Hol_{\Sigma}(\mathcal F(X),S)$ the holonomy of $\mathcal F(X)$
with respect to $S$ evaluated at a section $\Sigma$ transversal to
$S$. Notice that $\Sigma\cap S=\{p\}$ is a single point, $\Sigma$
is biholomorphic to a disc in ${\mathbb{C}}^{n-1}$ with center
corresponding to $p$ and $\Hol_{\Sigma
}({\mathcal F(X)},S)$ is conjugate to a subgroup of the group $\Diff({\mathbb{C}%
}^{n-1},0)$ of germs of complex diffeomorphisms fixing the origin
in ${\mathbb{C}}^{n-1}$. Further we  say that $f$ is
\textit{adapted} to $(\mathcal F(X),S)$ if it can be written
locally in the form $f=g/h$ where $g,h\in\mathcal{O}_{n}$ are
relatively prime, $S\subset Z(g)\cap Z(h)$, where $Z(g)$ and
$Z(h)$ denote the zero sets of $g$ and $h$ respectively, and
$\left. f\right\vert _{_{\Sigma}}$ is pure meromorphic for a
generic section $\Sigma$ transversal to $S$. Denote the set of
non-negative integers by $\mathbb{N}$ and the set of positive
integers by $\mathbb{Z}_{+}$. Let
$x=(x_{1},\cdots,x_{n})\in\mathbb{C}^{n}$ and
$I=(i_{1},\cdots,i_{n})\in\mathbb{N}^{n}$ then we  use
respectively
the following notation for monomials and their orders: $x^{I}:=x_{1}^{i_{1}%
}\cdots x_{n}^{i_{n}}$, and $|I|=i_{1}+\cdots+i_{n}$. Finally, we
denote by $C_{n}$ the germ of curve given by the union of the
coordinate axes of $(\mathbb{C}^{n},0)$. Given a vector field  $X$
in $\Gen(\frak X(\bc^3,0))$ and a nonvanishing holomorphic
function $u$ in a neighborhood of the origin $0\in\mathbb C^3$,
the vector field $Y=u X$ also satisfies $Y \in \Gen(\frak
X(\bc^3,0))$. We shall then say that $X$ and $Y$ are {\it
tangent}. As we will see (cf. Proposition \ref{first jet}), any
integrable vector field $X\in\Gen(\frak X(\bc^3,0))$ satisfies the
following condition:

\begin{Definition}
\label{definition:conditionstar}  {\rm  Let $X \in \Gen(\frak
X(\bc^3,0))$. We say that $X$ {\it satisfies condition} $(\star)$
if there is a real line $L\subset \mathbb R^2$ through the origin
containing the eigenvalues of $X$ such that one of the connected
components $L\setminus \{0\}$ contains a single eigenvalue
$\lambda(X)$ of $X$. In other words, not all the eigenvalues of
$X$ belong to the same connected component of $L\setminus \{0\}$.
The above condition holds for $X$ if and only if holds for any
vector field $Y$ such that $X$ and $Y$ are tangent. Condition
$(\star)$ implies that $X$ is in the Siegel domain, but is
stronger than this last. If $X$ satisfies $(\star)$
 we denote by $S_X$ the smooth invariant curve
associate to
 $\lambda(X)$.}
 \end{Definition}

 Our main result reads as follows:

\begin{Theorem}
\label{topological criterion}Suppose that $X\in\Gen(\frak
X(\bc^3,0))$ satisfies condition $(\star)$ and let $S_X$ be the
smooth invariant curve associate with the eigenvalue $\lambda(X)$.
Then the following conditions are equivalent:

\begin{enumerate}
\item[{\rm (i)}]  The leaves of $\mathcal F(X)$ are closed off $\operatorname{{sing}%
}(\mathcal F(X))$;

\item[{\rm (ii)}] $\left.  \Hol\right\vert _{\Sigma}(\mathcal
F(X),S_X)$ has finite orbits;

\item[{\rm (iii)}] $\left.  \Hol\right\vert _{\Sigma}(\mathcal
F(X),S_X)$ is periodic {\rm(}in particular linearizable and
finite{\rm)};

\item[{\rm (iv)}] $\mathcal F(X)$ has a holomorphic first
integral.
\end{enumerate}
\end{Theorem}

As a straightforward consequence we obtain
the following topological criterion for the existence of $\mathcal F(X)%
$-invariant meromorphic functions in $\Gen(\frak X(\bc^3,0))$.

\begin{Theorem}
\label{F-invariant merom.} Let $X\in\Gen(\frak X(\bc^3,0))$ and
$S_X$ be the invariant smooth curve associated with the eigenvalue
$\lambda(X)$. Suppose that $\mathcal F(X)$ has closed leaves off
$\operatorname{{sing}}({\mathcal F(X)})$, then there is an $\mathcal F(X)%
$-invariant meromorphic function adapted to $(\mathcal F(X),S_X)$.
\end{Theorem}

\section{Holomorphic first integrals and first jets}

Here we study the necessary conditions on the eigenvalues of a
vector field $X$ in $\Gen(\frak X(\bc^3,0))$ in order to have $X$
integrable. In particular we determine a criterion for linear
vector fields to have holomorphic first integrals.

\subsection{Algebraic criterion}

Let us give an algebraic description of the first jets of
integrable vector fields. First recall that $X\in\Gen(\frak
X(\bc^3,0))$ has a holomorphic first integral $F=(f_{1},f_{2})$ if
and only $i_{X}df_{j}=0$, $j=1,2$, and $f_{1},f_{2}$ are
transversal off the singular set of $X$. Let $X$ be given in the
system of coordinates $(x_{1},x_{2},x_{3})$ of
$(\mathbb{C}^{3},0)$ by%
\[
X(x)=\lambda_{1}x_{1}(1+a_{1}(x))\frac{\partial}{\partial x_{1}}+\lambda
_{2}x_{2}(1+a_{2}(x))\frac{\partial}{\partial x_{2}}+\lambda_{3}x_{3}%
(1+a_{3}(x))\frac{\partial}{\partial x_{3}}%
\]
where $a_{1},a_{2},a_{3}\in\mathcal{M}_{3}$ and consider a
holomorphic function $f\in\mathcal{M}_{3}$ preserving the
coordinate axes. Then $f$ must be given by $f(x)=\sum_{\left\vert
N\right\vert \geq p}a_{N}x^{N}$  where $a_{N}\in\mathbb{C}$,
$p\geq2$ and $N\in\mathbb{N}^{3}-C_{3}$. Since
$\frac{\partial f}{\partial x_{j}}=\sum_{|N|\geq p}n_{j}a_{N}x^{N}x_{j}^{-1}$, then%

\begin{align*}
J^{p}(df(X))  &  =\frac{\partial f(x)}{\partial x_{1}}\cdot(\lambda_{1}%
x_{1})+\frac{\partial f(x)}{\partial x_{2}}\cdot(\lambda_{2}x_{2}%
)+\frac{\partial f(x)}{\partial x_{3}}\cdot(\lambda_{3}x_{3})\\
&  =\sum_{\left\vert N\right\vert \geq p}(\lambda_{1}n_{1}+\lambda_{2}%
n_{2}+\lambda_{3}n_{3})a_{N}x^{N}\text{.}%
\end{align*}

Hence, if $i_{X}df=0$ then%
\begin{equation}
0=(\lambda_{1}n_{1}+\lambda_{2}n_{2}+\lambda_{3}n_{3})a_{N}\text{ for all
}|N|=p\text{, }N\in\mathbb{N}^{3}-C_{3}\text{. } \label{eq8}%
\end{equation}

From (\ref{eq8}) it follows that $a_{N}=0$ whenever $\lambda_{1}n_{1}%
+\lambda_{2}n_{2}+\lambda_{3}n_{3}\neq0$, for all $N\in\mathbb{N}^{3}-C_{3}$
such that $|N|=p$. In particular, in the absence of a resonance of the form%
\begin{equation}
\ \lambda_{1}n_{1}+\lambda_{2}n_{2}+\lambda_{3}n_{3}=0\text{,}
\label{resonance}%
\end{equation}
there will be no first integrals. Thus we have to study under what
conditions on the eigenvalues of $X$  we have a holomorphic first
integral of the form $F=(f_{1},f_{2})$ where $\ker(df_{1})$ is
transverse to $\ker(df_{2})$ off the singular set of $X$. For
further purposes we prove the following technical result.

\begin{Lemma}
\label{preparation}Let $\lambda_{1},\lambda_{2},\lambda_{3}\in\mathbb{C}%
^{\ast}$, and let $(n_{1},n_{2},n_{3}),(m_{1},m_{2},m_{3})\in\mathbb{N}%
^{3}-C_{3}$ be linearly independent and satisfying \textrm{(\ref{resonance})}
above. Then there are $m,n,k\in\mathbb{Z}$ and $\lambda\in\mathbb{C}^{\ast}$
such that
\[
(\lambda_{1},\lambda_{2},\lambda_{3})=\lambda(m,n,k)
\]
and $m\cdot n\cdot k<0$.
\end{Lemma}

\begin{proof}
Since $(n_{1},n_{2},n_{3}),(m_{1},m_{2},m_{3})$ are linearly independent then
$(n_{1}m_{2}-n_{2}m_{1},n_{1}m_{3}-n_{3}m_{1},n_{2}m_{3}-n_{3}m_{2}%
)\neq(0,0,0)$. We may suppose without loss of generality that $n_{1}%
m_{2}-n_{2}m_{1}\neq0$. Now consider the system of equations \
\begin{equation}
\left\{
\begin{array}
[c]{c}%
n_{1}\lambda_{1}+n_{2}\lambda_{2}+n_{3}\lambda_{3}=0\\
m_{1}\lambda_{1}+m_{2}\lambda_{2}+m_{3}\lambda_{3}=0
\end{array}
\right.  \label{eq1}%
\end{equation}
and sum the first equation multiplied by $-m_{1}$ with the second one
multiplied by $n_{1}$ to obtain $(n_{1}m_{2}-n_{2}m_{1})\lambda_{2}%
+(n_{1}m_{3}-n_{3}m_{1})\lambda_{3}=0$, i.e., $\lambda_{2}=-\left(
\frac{n_{1}m_{3}-n_{3}m_{1}}{n_{1}m_{2}-n_{2}m_{1}}\right)
\lambda_{3}$. Back to the first equation of (\ref{eq1}) we  have
$\lambda_{1}=\left(
\frac{n_{2}m_{3}-n_{3}m_{2}}{n_{1}m_{2}-n_{2}m_{1}}\right)
\lambda_{3}$. Hence
$(\lambda_{1},\lambda_{2},\lambda_{3})=\lambda(m,n,k)$ where
$\lambda=\frac{\lambda_{3}}{n_{1}m_{2}-n_{2}m_{1}}$, $m=n_{2}m_{3}-n_{3}m_{2}%
$, $n=n_{3}m_{1}-n_{1}m_{3}$ and $k=n_{1}m_{2}-n_{2}m_{1}$. Thus%
\[
mnk=(m_{1}m_{2}m_{3})^{2}(\frac{n_{2}}{m_{2}}-\frac{n_{3}}{m_{3}})(\frac
{n_{3}}{m_{3}}-\frac{n_{1}}{m1})(\frac{n_{1}}{m_{1}}-\frac{n_{2}}{m_{2}%
})\text{.}%
\]
An elementary evaluation shows that if two of the factors in the product
$(\frac{n_{2}}{m_{2}}-\frac{n_{3}}{m_{3}})(\frac{n_{3}}{m_{3}}-\frac{n_{1}%
}{m1})(\frac{n_{1}}{m_{1}}-\frac{n_{2}}{m_{2}})$ are positive, then the third
one must be negative.
\end{proof}

\begin{Proposition}
\label{first jet}Suppose that $X\in\Gen(\frak X(\bc^3,0))$ has a
holomorphic first integral, then $\mathcal F(X)$ can be given in
local coordinates by a vector field of the form
\[
X(x)=mx_{1}(1+a_{1}(x))\frac{\partial}{\partial x_{1}}+nx_{2}(1+a_{2}%
(x))\frac{\partial}{\partial x_{2}}-kx_{3}(1+a_{3}(x))\frac{\partial}{\partial
x_{3}}%
\]
where $m,n,k\in\mathbb{Z}_{+}$ and
$a_{1},a_{2},a_{3}\in\mathcal{M}_{3}$. In particular $X$ satisfies
condition $(\star)$.
\end{Proposition}

\begin{proof}
Suppose that $J^{1}(X)=\lambda_{1}x_{1}\frac{\partial}{\partial x_{1}}%
+\lambda_{2}x_{2}\frac{\partial}{\partial x_{2}}+\lambda_{3}x_{3}%
\frac{\partial}{\partial x_{3}},$ then Lemma \ref{preparation} assures that
its enough to prove that there is a pair of linearly independent vectors
$M,N\in\mathbb{N}^{3}-C_{3}$ satisfying (\ref{resonance}). Suppose that
$F=(f,g)$ is a first integral for $X$, with $f(x)=\sum_{\left\vert
N\right\vert \geq p}a_{N}x^{N}$ and $g(x)=\sum_{\left\vert N\right\vert \geq
q}b_{N}x^{N}$, i.e. with orders respectively $p,q$ ($\geq2$). From
(\ref{resonance}) we have that $0=(\lambda_{1}n_{1}+\lambda_{2}n_{2}%
+\lambda_{3}n_{3})a_{N}$ for all $|N|=p$. If there are two distinct
$a_{N},a_{N^{\prime}}\neq0$ then $N$ and $N^{\prime}$ satisfy the desired
condition. Reasoning in the same manner for $g$ we just have to consider the
case where $p=q$ and the terms of degree $p$ of $f$ and $g$ are monomials,
i.e., $f(x)=a_{P}x^{P}+\sum_{\left\vert N\right\vert \geq p+1}a_{N}x^{N}$ and
$g(x)=b_{P}x^{P}+\sum_{\left\vert N\right\vert \geq p+1}b_{N}x^{N}$ with
$|P|=p$, and $a_{P},b_{P}\neq0$. Now let $f_{1}:=\frac{1}{a_{P}}f-\frac
{1}{b_{P}}g$, then it can be written in the form $f_{1}(x)=h_{1}(x^{P}%
)+\sum_{\left\vert N\right\vert =q,N\notin\left\langle P\right\rangle }%
c_{N}x^{N}+h.o.t$, where $h\in\mathcal{M}_{1}$ is a polynomial
such that $\tau_{1}:=\operatorname{{deg}}(h_{1})<q$, where $q=|N|$
is the less natural number such that there exists $c_{N}\neq0$
form some $N\notin\left\langle P\right\rangle $ (notice that such
$q$  exists, since $f$ and $g$ are
transversal off the origin). Pick inductively $f_{k}:=f_{k-1}-h_{k-1}%
^{(\tau_{k-1})}(0)\left(  \frac{1}{b_{P}}g\right)  ^{\tau_{k-1}}$,
where $\tau_{k}:=\operatorname{{deg}}(f_{k})$, then after
repeating this process a finite number of steps we  have
$k_{0}\in\mathbb{Z}_{+}$ such that $f_{k_{0}}(x)=\sum_{\left\vert
N\right\vert =q,N\notin\left\langle P\right\rangle
}c_{N}x^{N}+h.o.t$. Since the set of $\mathcal F(X)$-invariant
holomorphic functions is a sub-ring of $\mathcal{O}_{3}$, then
$f_{k_{0}}$ is an $\mathcal F(X)$-invariant holomorphic function;
in particular it  satisfies (\ref{eq8}). Therefore, it is enough
to pick $R\notin\left\langle P\right\rangle $ such that $|R|=q$
and $c_{R}\neq0$.
\end{proof}

\subsection{Linear vector fields}

Let $X\in\Gen(\frak X(\bc^3,0))$ be a linear vector field, then we
determine under what conditions such a vector field has a
holomorphic first integral.  Recall from Proposition \ref{first
jet} that we just have to consider linear
vector fields of the form $X(x)=mx_{1}\frac{\partial}{\partial x_{1}}%
+nx_{2}\frac{\partial}{\partial
x_{2}}-kx_{3}\frac{\partial}{\partial x_{3}}$ where
$m,n,k\in\mathbb{N}$, and $m+n+k\geq2$. Thus we have to study
under what conditions on $m,n,k$ we  have a holomorphic first
integral $F$ for $X$. Write $F=(f_{1},f_{2})$ where $\ker(df_{1})$
is transverse to $\ker(df_{2})$ off the singular set of $X$. It is
well-known  that a germ of holomorphic foliation $\mathcal F(X)$
with $\operatorname*{codim}(\sing(\mathcal F(X)))=1$ can be
extended to a foliation $\operatorname*{Sat}(\mathcal F(X))$
called the \textit{saturated} of $\mathcal F(X)$ such that
$\codim(\sing(\operatorname*{Sat}(\mathcal{F)}))\geq2$.

\begin{Lemma}
\label{transversal}Let $N,M\in\mathbb{N}^{3}-C_{3}$ be two vectors
satisfying {\rm(\ref{resonance})}, and let $f(x)=x^{N}$,
$g(x)=x^{M}$. Then $\operatorname*{Sat}(df=0)$ is transversal to
$\operatorname*{Sat}(dg=0)$ if and only if $N$ and $M$ are
linearly independent.
\end{Lemma}

\begin{proof}
First suppose that the vectors are linearly dependent, then there is
$k\in\mathbb{N}$ such that $m_{j}=kn_{j}$ for all $j=1,2,3$. In particular
$g=(f)^{k}$ and thus $dg=(f)^{k-1}df$. Therefore whenever $\operatorname*{Sat}%
(dg=0)$ is non-singular it coincides with
$\operatorname*{Sat}(df=0)$. Conversely suppose that the vectors
are linearly independent, then we have to show that
$\operatorname*{Sat}(df=0)$ is transverse to
$\operatorname*{Sat}(dg=0)$ away from the origin, or equivalently
that $\operatorname*{Sat}(df\wedge dg)$ in
non-singular off the origin. In fact,%
\begin{align*}
df\wedge dg  &  =x_{1}^{n_{1}+m_{1}-1}x_{2}^{n_{2}+m_{2}-1}x_{3}^{n_{3}%
+m_{3}-1}[(n_{1}m_{2}-n_{2}m_{1})x_{3}dx_{1}\wedge dx_{2}\\
&  +(n_{1}m_{3}-n_{3}m_{1})x_{2}dx_{1}\wedge dx_{2}+(n_{2}m_{3}-n_{3}%
m_{2})x_{1}dx_{2}\wedge dx_{3}]
\end{align*}
Hence $\Sat(df\wedge dg)$ is given by $\omega=0$ where $\omega=(n_{1}%
m_{2}-n_{2}m_{1})x_{3}dx_{1}\wedge dx_{2}+(n_{1}m_{3}-n_{3}m_{1})x_{2}%
dx_{1}\wedge dx_{2}+(n_{2}m_{3}-n_{3}m_{2})x_{1}dx_{2}\wedge dx_{3}$ vanishes
at the origin and whenever $n_{1}m_{2}-n_{2}m_{1}=n_{1}m_{3}-n_{3}m_{1}%
=n_{2}m_{3}-n_{3}m_{2}=0$. But the later happens exactly when $(n_{1}%
,n_{2},n_{3})$ and $(m_{1},m_{2},m_{3})$ are linearly dependent.
\end{proof}

An algebraic characterization of integrable linear vector fields
is given by the following result.

\begin{Lemma}
\label{1st integral}Any linear vector field of the form $X(x)=mx_{1}%
\frac{\partial}{\partial x_{1}}+nx_{2}\frac{\partial}{\partial x_{2}}%
-kx_{3}\frac{\partial}{\partial x_{3}}$, where $(m,n,k)\in\mathbb{Z}_{+}^{3}$,
has a holomorphic first integral of the form $F(x)=(x^{N},\allowbreak x^{M})$,
where $N,M\in\mathbb{N}^{3}-C_{3}$.
\end{Lemma}

\begin{proof}
From Lemma \ref{transversal} and the calculations made in order to
obtain (\ref{eq8}), one can easily check that this is just a
matter of finding two linearly independent solutions in
$\mathbb{N}^{3}-C_{3}$ for the homogeneous equation $\
mx+ny-kz=0$. Therefore, we just have to pick $x_{j}:=k\widetilde
{x}_{j}$ and $y_{j}:=k\widetilde{y}_{j}$, $j=1,2$, where $(\widetilde{x}%
_{1},\widetilde{y}_{1}),(\widetilde{x}_{2},\widetilde{y}_{2})\in\mathbb{N}%
^{2}$ are linearly independent.
\end{proof}

\subsection{Finite subgroups of germs of diffeomorphisms}

Let $f\colon U\subset \bc^n \to f(U)\subset \bc^n$ be a
homeomorphism in a neighborhood of the origin $0\in\bc^n$ with
$f(0) = 0$. Given  $x\in U$ denote by $\O_U(f, x)$ the $f$-orbit
of $x$ that does not leave $U$, i.e. $y\in \O_U(f, x)$ if and only
if $\{x, f(x), . . . , y = f[k](x)\}\subset U$ or $\{x,
f^{[-1]}(x), . . . , y = f^{[-k]}(x)\}\subset U$ for some $k\in
\mathbb N$. The following result is very useful (cf. Theorem 3.1,
\cite{Ma 2003}):

\begin{Theorem}
\label{Theorem:brochero} Let $F\in \Diff(\bc^2,0)$. The group
generated by $F$ is finite if and only if there exists a
neighborhood $V$ of $0\in \bc^2$, such that: \begin{enumerate}

\item $|O_V (F,X)| < \infty$ for all $X\in V$ and

\item $F$ leaves invariant infinitely many analytic varieties at
$0$.

\end{enumerate}

\end{Theorem}

\section{Holomorphic first integrals
 in $\Gen(\frak X(\bc^3,0))$}

In this section we  describe a topological criterion for the
integrability of $\Gen(\frak X(\bc^3,0))$ foliations.

\subsection{The transverse geometry of $\mathcal{F}$}

As a first step we  study the existence of a transversely
dicritical codimension $1$ foliation $\mathcal{G}$ tangent to the
orbits of a vector filed $X\in\Gen(\frak X(\bc^3,0))$.

\begin{definition}
{\rm Let $\mathcal{F}$ be a germ of foliation by curves on
$(\mathbb{C}^{3},0)$ given by the vector field $X$, then a germ of
codimension one foliation $\mathcal{G}$ on $(\mathbb{C}^{3},0)$ is
said to be {\it tangent} to $\mathcal{F}$ if its leaves contain
the orbits of $X$.}
\end{definition}

\begin{definition}
{\rm Let $\mathcal{G}$ be a codimension one foliation on
$(\mathbb{C}^{3},0)$, and  $S$ a germ of curve through the origin,
invariant by $\mathcal{G}$. We say that $\mathcal{G}$ is {\it
transversely dicritical} with respect to $S$ if for any section
$\Sigma$ transversal to $S$ in $(\mathbb{C}^{3},0)$ the
restriction $\left.  \mathcal{G}\right\vert _{\Sigma}$ is
dicritical.}
\end{definition}

\begin{Lemma}
\label{transversal dicrit.}Let $\mathcal{F}$ be a germ of
foliation by curves on $(\mathbb{C}^{3},0)$, $S$ an invariant
curve of $\mathcal{F}$ through the origin and $\mathcal{G}$ a
codimension one foliation satisfying the following conditions:

\begin{enumerate}
\item[{\rm (i)}] $\mathcal{G}$ is tangent to $X$;

\item[{\rm (ii)}] There is a section $\Sigma$ transversal to $S$
such that $\left. \mathcal{G}\right\vert _{\Sigma}$ is dicritical.
\end{enumerate}

Then $\mathcal{G}$ is transversely dicritical with respect to $S.$
\end{Lemma}

\begin{proof}
Since the orbits of $X$ are contained in the leaves of
$\mathcal{G}$, then they are invariant by the flow of $X$.
Therefore if $\Sigma^{\prime}$ is another section transversal to
$S$ and $\phi:\Sigma\longrightarrow \Sigma^{\prime}$ is an element
of the holonomy pseudogroup of $X$ with respect to $S$, then it is
a diffeomorphism taking the leaves of $\left.
\mathcal{G}\right\vert _{\Sigma}$ onto the leaves of $\left.  \mathcal{G}%
\right\vert _{\Sigma^{\prime}}$.
\end{proof}

Our main concern here is the following: given a germ of foliation
by curves
$\mathcal{F}$ induced by a germ of vector field of the form%
\begin{equation}
X=mx(1+a(x,y,z))\frac{\partial}{\partial
x}+ny(1+b(x,y,z))\frac{\partial
}{\partial x}-kz(1+c(x,y,z)\frac{\partial}{\partial z} \label{eq4}%
\end{equation}
with $a,b,c\in\mathcal{M}_{3}$, find a codimension $1$ germ of
holomorphic foliation tangent to $X\ $which is transversely
dicritical with respect to $S$. Thus, in order to achieve our task
we have to find a solution of the following problem:

\begin{Problem}
\label{geometric integrability}Find a germ of holomorphic one form
$\omega(x,y,z)=P(x,y,z)dx+Q(x,y,z)dy+R(x,y,z)dz$ satisfying the
following conditions:

\begin{enumerate}
\item[{\rm (i)}] {\rm(}integrability{\rm)} $\omega\wedge
d\omega=0$;

\item[{\rm (ii)}] {\rm(}tangency{\rm)} $i_{X}\omega=0$;

\item[{\rm (iii)}] {\rm(}dicriticalness{\rm)}
$(\mathcal{G}:\omega=0)$ is transversely dicritical with respect
to $S:=(z=0)$.
\end{enumerate}
\end{Problem}

\subsection{Perturbations of some dicritical foliations on $(\mathbb{C}%
^{2},0)$} \label{subsection:perturbations} In this section we
shall deal with the question of transversal dicriticity of a
codimension one foliation $\mathcal{G}:(\omega=0)$ with respect to
$S$. We begin with  a vector field $X\in\Gen(\frak X(\bc^3,0))$ of
the form $X=mx(1+P(x,y,z))\frac{\partial}{\partial
x}+ny(1+Q(x,y,z))\frac{\partial }{\partial
x}-kz(1+R(x,y,z)\frac{\partial}{\partial z}$ where
$m,n,k\in\mathbb N$ and $P, Q, R\in$ $\mathcal{M}_{3}$ we can
define the one-parameter family of $1$-forms
\begin{equation}
\omega_{z}(x,y):=\omega(x,y,z)=mx(1+P_{z}(x,y))dy-ny(1+Q_{z}%
(x,y))dx\label{eq4.1}%
\end{equation}
where $\omega_{z}\in\Lambda^{1}(\mathbb{C}^{2},0)$,
$P_{z}(x,y):=P_{z}(x,y,z)$ and $Q_{z}%
(x,y):=Q(x,y,z)$, for each fixed $z\in(\mathbb{C},0)$. We  denote
by $\mathcal F(X)_{z}^{\prime}$ the one-parameter family of germs
of foliations defined on $(\mathbb{C}^{2},0)$ by $\omega_{z}$.

\begin{Remark}
\label{tangency} {\rm Notice that for $X$ as above we have
$\lambda(X)=-k$ and $S_X$ is the $z$-axis. Also the kernel of
$\omega$ denoted  $\operatorname*{Ker}(\omega)$ is a holomorphic
distribution on $(\bc^3,0)$ tangent to the orbits of $X$, since
$i_{X}\omega=0$. Moreover,  if $\Sigma_{z}:=(z=const.)$ and
$i_{z}:\Sigma_{z}\longrightarrow\mathbb{C}^{3}$ is the inclusion
map, then $(i_{z})^{\ast}X$ is tangent to $(i_{z})^{\ast
}\omega$.}
\end{Remark}

Keeping in mind this framework,  we consider the following
situation.
Let $\omega_{z}$ be an one-parameter holomorphic family of $1$-forms on $(\mathbb{C}%
^{2},0)$ given by
\[
\omega_{z}(x,y):=mx(1+P_{z}(x,y))dy-ny(1+Q_{z}(x,y))dx
\]
where $P,Q\in\mathcal{M}_{3}$, $P_{z}(x,y):=P_{z}(x,y,z)$ and $Q_{z}%
(x,y):=Q(x,y,z)$, for each fixed $z\in(\mathbb{C},0)$. We shall
denote by $\mathcal{F}_{z}$ the one-parameter family of germs of
foliations defined on $(\mathbb{C}^{2},0)$ by $\omega_{z}$.

\begin{Lemma}
\label{point dicrit.} For each neighborhood $U$ of the origin in
$\mathbb{C}$, there is $z\in U$ such that $\mathcal{F}_{z}$ is
dicritical.
\end{Lemma}

\begin{proof}
First notice that $\mathcal{F}_{z}$ may be given by $\omega_{z}%
(x,y):=x(1+xA+yB)dy-y(\frac{p}{q}+g(z)+xC+yD)dx$, where $p$ and
$q$ are relatively prime positive integers, $A,B,C,D\in$
$\mathcal{O}_{3}$, and $g\in\mathcal{M}_{1}$. In particular, we
may suppose that $\left\vert g(z)\right\vert <\varepsilon$ for any
$z$ sufficiently close to the origin. Moreover, since $g$ is an
open map, we can find a sequence $z_{n}\rightarrow0$ such that
$g(z_{n})=\frac{1}{q^{n}}$. Therefore, we just need to verify that
$(\mathcal{F}_{z_{n}}:\omega_{z_{n}}=0)$ is dicritical; but this
comes immediately from Lemma \ref{dicritical resolution} below.
\end{proof}

\begin{Lemma}
\label{dicritical resolution}Let $\mathcal F(\omega)$ be a
foliation in $(\bc^2,0)$  be given by a one-form
$\omega(x,y)=mx(1+xA(x,y)+yB(x,y))dy-ny(1+xC(x,y)+yD(x,y))dx$
where $A,B,C,D\in\mathcal{O}_{2}$. Then $\mathcal F(\omega)$ has
the same resolution tree of $\mathcal F(\omega_{m,n})$ where
$\omega_{m,n}(x,y)=mxdy-nydx$. In particular $\mathcal F(\omega)$
is dicritical.
\end{Lemma}

\begin{proof}
After one blow-up we  obtain
\begin{equation}
\widetilde{\omega}_{0}(x,t)
=\frac{1}{x}\pi^{\ast}\omega(x,t)\nonumber
=x(m+x\widetilde{A}+tx^{2}\widetilde{B})dt+t((m-n)+x(\widetilde
{A}-\widetilde{C})+tx(\widetilde{B}-\widetilde{D}))dx \label{eq2}%
\end{equation}
and%
\begin{equation}
\widetilde{\omega}_{\infty}(u,y)    =\frac{1}{y}\pi^{\ast}\omega
(u,y)\nonumber\\
  =u((m-n)+uy(\widetilde{A}-\widetilde{C})+y(\widetilde{B}-\widetilde
{D}))dy-y(n+uy\widetilde{C}+y\widetilde{D})dt+ \label{eq3}%
\end{equation}

where $\widetilde{A}=\pi^{\ast}A,\widetilde{B}=\pi^{\ast}B,\widetilde{C}%
=\pi^{\ast}C$ and $\widetilde{D}=\pi^{\ast}D$ . Without loss of
generality we may suppose that $m>n$. In this case we have a
reduced singularity at $(t,x)=0$ and a non-reduced one at
$(u,y)=(0,0)$, with the same character of the original one.
Therefore, Euclid's algorithm assures that after a finite number
of blow-ups we  obtain a linear chain of projective lines, any of
them with two singularities, and just one of the singularities in
the exceptional divisor being non-reduced: the one given in local
coordinates $(\widetilde{x},\widetilde{y})$ by
$\omega_{m,n}(\widetilde{x},\widetilde
{y})=\widetilde{x}(1+\widetilde{x}A_{m,n}(\widetilde{x},\widetilde
{y})+\widetilde{y}B_{m,n}(\widetilde{x},\widetilde{y}))d\widetilde
{y}-\widetilde{y}(1+\widetilde{x}C_{m,n}(\widetilde{x},\widetilde
{y})+\widetilde{y}D_{m,n}(\widetilde{x},\widetilde{y}))d\widetilde{x}$.
Hence (\ref{eq2}) and (\ref{eq3}) show that one further blow-up
leads to a foliation transversal to the new projective line
appearing in the divisor; this is exactly the resolution tree of
$\omega_{m,n}$.
\end{proof}

\subsection{Solution to Problem~\ref{geometric integrability}}
\label{subsection:solution}  We verify that we can find (local)
solutions for the Problem \ref{geometric integrability}.

\subsubsection{Tangency and dicriticity conditions}

First notice that Lemmas \ref{transversal dicrit.} and \ref{point
dicrit.} suggest us to search $\omega$ of the form $\omega
=ny(1+q(x,y,z))dx-mx(1+p(x,y,z))dy+R(x,y,z)dz$. Thus, from the
tangency condition of Problem \ref{geometric integrability}, we
can determine $R$ in
terms of $p$ and $q$ as follows. From (\ref{eq4}) we shall have $0=i_{X}%
\omega=mnxy[a-b+q-p+aq-bp]-kz(1+c)R$. If we set $K:=a-b+q-p+aq-bp$
we shall
have%
\begin{equation}
0=mnxyK-kz(1+c)R\label{eq14}%
\end{equation}
Recall that our main goal here is to obtain $R$ as an explicit
holomorphic function in terms of the coordinates of $X$ and the
first two coordinates of $\omega$. Therefore, we shall need $K$\
to be a multiple of\ $z$. But notice that, in the particular case
where $C\equiv0$ we just have to take $p=a$ and $q=b$. This
suggests us to take $q:=b+z\widetilde{q}$ and $p:=a+z\widetilde
{p}$, for some $\widetilde{p},\widetilde{q}\in\mathcal{O}_{3}$.
Under this assumption we shall have after some simple calculations
that $K=z[(a+1)\widetilde{q}-(b+1)\widetilde{p}]$. Therefore
(\ref{eq14}) turns out
to be $0=mnxyz[(a+1)\widetilde{q}-(b+1)\widetilde{p}]-kz(1+c)R$. Thus%
\begin{equation}
R=\frac{mnxy[(a+1)\widetilde{q}-(b+1)\widetilde{p}]}{k(1+c)}\text{.}%
\label{eq15}%
\end{equation}

\subsubsection{Integrability condition}

First notice that for $\omega=P(x,y,z)dx+Q(x,y,z)dy+R(x,y,z)dz$
the integrability condition $\omega\wedge d\omega=0$ is equivalent
to the first order PDE
\begin{equation}
(-P_{y}+Q_{x})R-(-P_{z}+R_{x})Q+(-Q_{z}+R_{y})P=0,\label{eq12}%
\end{equation}
therefore we have to apply (\ref{eq15}) in (\ref{eq12}) in order
to study the solutions of Problem \ref{geometric integrability}.
For simplicity, let us first introduce the following notation:
$\overline{a}:=1+a,\overline
{b}:=1+b,\overline{c}:=1+c,\widetilde{q}:=\overline{c}\overline{q}$
and $\widetilde{p}:=\overline{c}\overline{p}$. Then, if we take
$P,Q$ and $R$ in
the form%
\begin{equation}
\left\{
\begin{array}
[c]{l}%
P=ny(\overline{b}+z\overline{c}\overline{q})\text{,}\\
Q=-mx(\overline{a}+z\overline{c}\overline{p})\text{,}\\
R=\frac{mn}{k}xy(\overline{a}\overline{q}-\overline{b}\overline{p})\text{,}%
\end{array}
\right.  \label{eq16}%
\end{equation}
we shall obtain after some lengthly calculations that%
\begin{align*}
0 &  =mx[(\overline{a}^{2}\overline{q}_{x}-\overline{a}\overline{b}%
_{x}\overline{p}-\overline{a}\overline{b}\overline{p}_{x}+z\overline{a}%
_{x}\overline{c}\overline{p}\overline{q}+z\overline{a}\overline{c}\overline
{p}\overline{q}_{x}-z\overline{b}_{x}\overline{c}\overline{p}^{2})\\
&
-(z\overline{a}\overline{c}_{x}\overline{p}\overline{q}+z\overline
{a}\overline{c}\overline{p}_{x}\overline{q}-\overline{a}_{x}\overline
{b}\overline{p}-z\overline{b}\overline{c}_{x}\overline{p}^{2})]\\
&
+ny[(\overline{b}\overline{a}_{y}\overline{q}+\overline{b}\overline
{a}\overline{q}_{y}-\overline{b}^{2}\overline{p}_{y}+z\overline{a}%
_{y}\overline{c}\overline{q}^{2}-z\overline{b}_{y}\overline{c}\overline
{p}\overline{q}-z\overline{b}\overline{c}\overline{q}\overline{p}_{y})\\
&  -(\overline{a}\overline{b}_{y}\overline{q}+z\overline{a}\overline{c}%
_{y}\overline{q}^{2}-z\overline{b}\overline{c}_{y}\overline{p}\overline
{q}-z\overline{b}\overline{c}\overline{p}\overline{q}_{y})]\\
&
+k[(\overline{a}_{z}\overline{b}+\overline{b}\overline{c}\overline
{p}+z\overline{b}\overline{c}_{z}\overline{p}+z\overline{b}\overline
{c}\overline{p}_{z}+z\overline{a}_{z}\overline{c}\overline{q}+z^{2}%
\overline{c}^{2}\overline{p}_{z}\overline{q})\\
&
-(\overline{a}\overline{b}_{z}+\overline{a}\overline{c}\overline
{q}+z\overline{a}\overline{c}_{z}\overline{q}+z\overline{a}\overline
{c}\overline{q}_{z}+z\overline{b}_{z}\overline{c}\overline{p}+z^{2}%
\overline{c}^{2}\overline{p}\overline{q}_{z})]
\end{align*}
Since we are interested in the foliation defined by $\omega$ (and
not in $\omega$ itself) and since $\overline{b}$ is a unity, then
we may suppose without loss of generality that $\overline{q}=0$
(In fact, from (\ref{eq16}) we just have to divide $\omega$ by
$(1+z\overline{c}\overline{q}/\overline {b})$). Analogously, we
are not interested in $X$ itself but in the foliation
$\mathcal{F}$ defined by it. Thus we may also suppose without loss
of generality that $\overline{b}=1$ (just divide $X$ by
$\overline{b}=1+b$). Under these assumptions the above equation
turns out to be
\begin{equation}
0=ka_{z}+(mx\overline{a}_{x}+k\overline{c}+kz\overline{c}_{z})\overline
{p}+mxz\overline{c}_{x}\overline{p}^{2}-mxa\overline{p}_{x}-ny\overline{p}%
_{y}+kz\overline{c}\overline{p}_{z}\label{eq5}%
\end{equation}

Now notice that a germ of holomorphic function defined on any of
the coordinate planes can be taken as a Cauchy condition for the
first order PDE (\ref{eq5}). Since the origin is a
non-characteristic point for (\ref{eq5}), then we may apply the
classical method of characteristics for non-linear PDEs in order
to find a (local) solution to (\ref{eq5}) with the above initial
condition (cf. e.g. \cite{John 75}). This solves Problem
\ref{geometric integrability}, as desired.

\subsection{The existence of a topological criterion}

Here we prove Theorem \ref{topological criterion}. For this sake
let us first recall some facts proved along this article and
introduce some terminology. Let $X \in\Gen(\frak X(\bc^3,0))$
satisfying condition $(\star)$ in \S~\ref{subsection:notation}. We
can assume that the curve $S_X$ is the $z$-axis so that we can
apply the results in \S~\ref{subsection:perturbations}. Let
$\Sigma_{z}:=(z=const.)$
be a section transversal to $S_X$ and $\Hol_{\Sigma_{z}%
}(\mathcal F(X),S_X)$  the holonomy of $\mathcal F(X)$ with
respect to $S_X$, evaluated at $\Sigma_{z}$.

\begin{proof}
[Proof of Theorem \ref{topological criterion}]First let us prove
that (i) implies (ii). Suppose all the leaves of $\mathcal F(X)$
are closed off $\sing({\mathcal F(X)})$. Given a leaf $L$ of
$\mathcal F$ it follows that the closure $\ov{L}\subset L\cup
\sing(\mathcal F)$ is an analytic subset of pure dimension one in
$\bc^3$ and since this leaf is transversal  to $\Sigma_{z}$, we
conclude that $\ov{L}\cap\Sigma_z$ is a finite set. On the other
hand, given a point $x\in L\cap \Sigma_z$ its orbit in the
holonomy group is also contained in $L\cap \Sigma_z$ so that it is
a finite set. Thus the  orbits of $H_{z}$ are finite proving (ii).
Now let us verify that that (ii) implies (iii). First recall from
Lemma \ref{point dicrit.} that there is $z_{0}\in(\mathbb{C},0)$
may be choose conveniently in such a way that $\mathcal
F(X)_{z_{0}}^{\prime}$ is dicritical. Now consider a simple loop
$\gamma$ around the origin, inside the $z$-axis, starting from
$z_{0}$. Pick a leaf $L$ of $\mathcal F(X)_{z_{0}}^{\prime}$ and
consider the liftings of $\gamma$ starting at points of $L$, along
the trajectories of $\mathcal F(X)$. Then these liftings form a
three dimensional real manifold, say $S_{L}$, whose intersection
with $\Sigma_{z_{0}}$ is given by $L$ and $L^{\prime}$ (see Figure
$1$). In particular if $h:=h_{\gamma}$ is the generator of
$\Hol_{\Sigma_{z}}(\mathcal F(X),S_X)$, then $L^{\prime }=h(L)$.
If $\omega$ is a solution to Problem \ref{geometric
integrability}, then $S_{L}$ is tangent to $\Ker(\omega)$, and
$S_{L}\cap$ $\Sigma_{z_{0}}$ in tangent to the distribution
$\Ker(\omega_{z_{0}})$.
 But the last one is integrable
and given by $\mathcal F(X)_{z_{0}}^{\prime}$, then $L^{\prime}$
is a leaf of $\mathcal F(X)_{z_{0}}^{\prime}$. Thus, since
$\mathcal F(X)_{z_{0}}^{\prime}$ is dicritical and the orbits of
$h$ are finite, then $h$ preserves an infinite number of
separatrices in $\Sigma_{z_{0}}$. Therefore,
Theorem~\ref{Theorem:brochero} assures that $h$ is periodic (in
particular linearizable and finite). Now let us prove that (iii)
implies (iv). From \cite{EliYa 84} this comes immediately from
Proposition \ref{first jet} and the fact that $\mathcal F(X)\in$
$\Gen(\frak X(\bc^3,0))$, since the holonomy of $\mathcal F(X)$ is
linearizable. Therefore one may suppose without loss of generality
that $X(x)=mx_{1}\frac{\partial}{\partial
x_{1}}+nx_{2}\frac{\partial}{\partial
x_{2}}-kx_{3}\frac{\partial}{\partial x_{3}}$ and the result
follows from Lemma \ref{1st integral}. Finally, since analytic
varieties are closed, then (i) comes straightforward from (iv).
\end{proof}

%

%TCIMACRO{\FRAME{dtbpFU}{3.039in}{1.9406in}{0pt}{\Qcb{Figure 1. The liftings of
%$\gamma$ along the leaves of $\QTR{cal}{F}$ starting at points of $L$.}}%
%{}{sl.eps}{\special{ language "Scientific Word";  type "GRAPHIC";
%width 3.039in;  height 1.9406in;  depth 0pt;  original-width 11.3057in;
%original-height 7.1805in;  cropleft "0";  croptop "1";  cropright "1";
%cropbottom "0";  filename 'SL.eps';file-properties "XNPEU";}}}%
%BeginExpansion
\begin{center}
\includegraphics[height=1.9406in,width=3.039in]{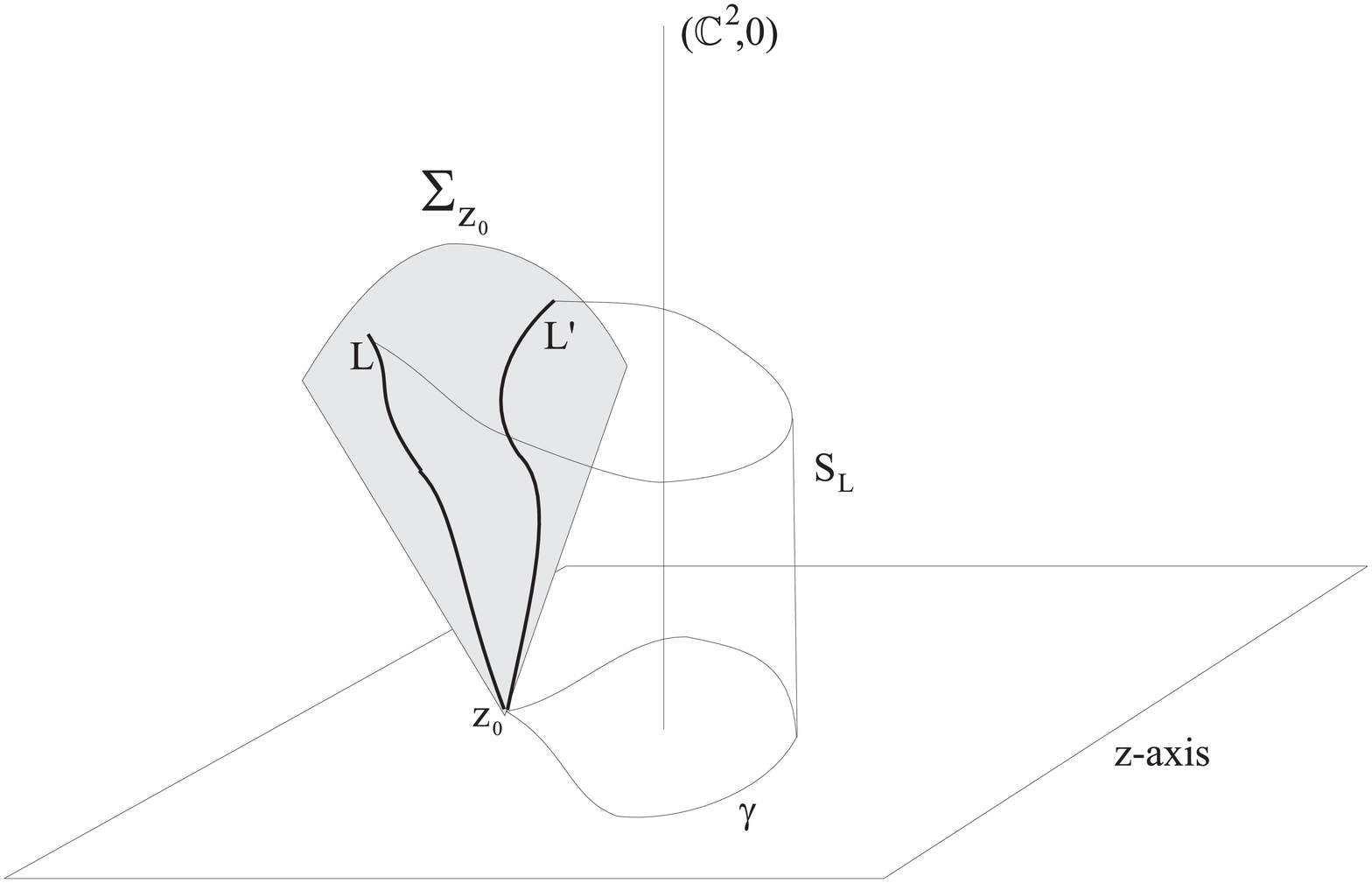}

Figure 1. The liftings of $\gamma$ along the leaves of $\mathcal
F(X)$ starting at points of $L$.
\end{center}
%EndExpansion

\vspace{0.1in}

\begin{proof}
[Proof of Theorem \ref{F-invariant merom.}] Suppose that
$X\in\Gen(\frak X(\bc^3,0))$ has all its leaves closed off the
singular set, then by Theorem \ref{topological criterion} we may
assume without loss of generality that $\
X(x)=mx_{1}\frac{\partial}{\partial x_{1}}+nx_{2}\frac{\partial
}{\partial x_{2}}-kx_{3}\frac{\partial}{\partial x_{3}}$ and has a
holomorphic first integral of the form
$F(x_{1},x_{2},x_{3})=(x^{P},x^{Q})$ for some linearly independent
$P,Q\in\mathbb{N}^{3}-C_{3}$ with $|P|,|Q|\geq2$. In particular we
have
\[
\left\{
\begin{array}
[c]{c}%
p_{1}m+p_{2}n-p_{3}k=0\\
q_{1}m+q_{2}n-q_{3}k=0
\end{array}
\right.
\]
and thus $(p_{1}q_{3}-p_{3}q_{1})m+(p_{2}q_{3}-p_{3}q_{2})n=0$. Since $m,n>0$
and $P,Q$ are linearly independent, then $(p_{1}q_{3}-p_{3}q_{1})(p_{2}%
q_{3}-p_{3}q_{2})<0$ . Therefore $f(x):=x^{q_{3}\cdot P}/x^{-p_{3}\cdot Q}$ is
an $\mathcal F(X)$-invariant meromorphic function adapted to $(\mathcal F(X)%
,S_X)$.
\end{proof}

\begin{tabular}
[c]{ll}%
Leonardo C\^amara & \qquad Bruno Sc\'ardua\\
Departamento de Matem\'atica - CCE & \qquad Instituto de Matem\'atica\\
Universidade Federal do Esp\'irito Santo & \qquad Universidade Federal do Rio
de Janeiro\\
Av. Fernando Ferrari 514 & \qquad Caixa Postal 68530\\
29075-910 - Vit\'oria - ES & \qquad21.945-970 Rio de Janeiro-RJ\\
BRAZIL & \qquad BRAZIL\\
camara@cce.ufes.br & \ scardua@impa.br
\end{tabular}

\end{document}